\documentclass[preprint,12pt]{elsarticle}

\usepackage{amssymb}
\usepackage{amsbsy}
\usepackage{amsmath}
\usepackage{color}

\journal{\textsf{arXiv}, typeset with elsarticle.cls}

    \newcommand{\floor}[1]{\lfloor#1\rfloor}

    \newcommand{\EE}{\mathbb{E}}

    \newcommand{\x}{\boldsymbol{x}}

    \renewcommand{\Pr}{\operatorname{P}}

    \newcommand{\dto}{\xrightarrow{d}}
    
    \newcommand{\vto}{\xrightarrow{v}}

    \newcommand{\toi}{\to\infty}

    \newcommand{\rmd}{\mathrm{d}}

    \newcommand{\be}{\begin{equation}}
    \newcommand{\ee}{\end{equation}}
    
    \renewcommand{\le}{\leqslant}
    \renewcommand{\ge}{\geqslant}
    \renewcommand{\leq}{\le}

  \newtheorem{thm}{Theorem}[section]
  \newtheorem{lem}[thm]{Lemma}
  \newdefinition{rem}[thm]{Remark}
  \newproof{pf}{Proof}
  \newproof{pot1}{Proof of Lemma~\ref{l:first}}
  \newproof{pot2}{Proof of Theorem~\ref{t:FinMA}}
  \newproof{pot3}{Proof of Theorem~\ref{t:InfMA}}

\numberwithin{equation}{section}

\begin{document}

\begin{frontmatter}

\title{A limit theorem for moving averages in the $\alpha$--stable
domain of attraction}

\author[fn1]{Bojan Basrak}
\ead{bbasrak@math.hr}

\author[fn2]{Danijel Krizmani\'{c}}
\ead{dkrizmanic@math.uniri.hr}

\address[fn1]{University of Zagreb, Department of Mathematics, Bijeni\v{c}ka 30, 10000 Zagreb, Croatia}
\address[fn2]{University of Rijeka, Department of Mathematics, Radmile Matej\v{c}i\'{c} 2, 51000 Rijeka, Croatia}

\begin{abstract}

In the early 1990's, Avram and Taqqu showed that
regularly varying moving average processes
 with all coefficients nonnegative and the tail index
$\alpha$ strictly between 0 and 2 satisfy functional limit theorem.
They also conjectured
that an equivalent statement holds under a certain less restrictive assumption
on the coefficients, but in a
different topology on the space of c\'adl\'ag functions.
We give a proof of this result.
\end{abstract}


\end{frontmatter}

\section{Introduction}
\label{intro}

It is known that the partial sums of i.i.d. regularly varying sequences
with the tail index $\alpha \in (0,2)$
satisfy the functional limit theorem
with an $\alpha$--stable L\'evy process as a limit.
This was first shown by Skorohod in 1950's  using his concept of $J_1$ topology on the space
of c\'adl\'ag functions. For a nice contemporary presentation of
this result we refer to
Resnick~\cite{Resnick07}, Chapter 7.
Naturally, one wonders whether the same holds for other
stationary sequences.
It turns out that, by introducing
 other alternative topologies  on the same space, Skorohod
gave us  the right tools to study this issue.
This was first observed by Avram and Taqqu~\cite{AvTa92}
who showed that the functional limit theorem holds for regularly
varying moving average processes provided that they have all coefficients
of the same sign. They used Skorohod's $M_1$ topology to obtain
the result and made a further conjecture that a similar
theorem holds under a less restrictive assumption on the coefficients
of the moving average process, but in somewhat weaker $M_2$ topology.
The principle goal of our paper is to show that this is indeed true.
We start by stating the problem precisely.

In the sequel, $(Z_{i})_{i \in \mathbb{Z}}$ denotes an i.i.d. sequence of regularly varying random variables with index of regular variation $\alpha \in (0,2)$.
In particular, this means that
$$ \Pr(|Z_{i}| > x) = x^{-\alpha} L(x), \qquad x>0,$$
where $L$ is a slowly varying function at $\infty$.
Let $(a_{n})$ be a sequence of positive real numbers such that
\be\label{eq:niz}
n \Pr (|Z_{1}|>a_{n}) \to 1,
\ee
as $n \to \infty$. Regular
variation of $Z_{i}$ can be expressed in terms of
vague convergence of measures on $\EE = \overline{\mathbb{R}} \setminus \{0\}$: for $a_n$ as in
\eqref{eq:niz} and as $n \to \infty$,
\begin{equation}
  \label{eq:onedimregvar}
  n \Pr( a_n^{-1} Z_i \in \cdot \, ) \vto \mu( \, \cdot \,),
\end{equation}
with the measure $\mu$ on $\EE$  given by
\begin{equation}
\label{eq:mu}
  \mu(\rmd x) = \bigl( p \, 1_{(0, \infty)}(x) + r \, 1_{(-\infty, 0)}(x) \bigr) \, \alpha |x|^{-\alpha-1} \, \rmd x,
\end{equation}
where
\be\label{eq:pq}
p =   \lim_{x \to \infty} \frac{\Pr(Z_i > x)}{\Pr(|Z_i| > x)} \qquad \textrm{and} \qquad
  r =   \lim_{x \to \infty} \frac{\Pr(Z_i < -x)}{\Pr(|Z_i| > x)}.
\ee
Under these assumptions on the sequence $(Z_{i})$, we study
the moving average process of the form
$$ X_{i} = \sum_{j=-\infty}^{\infty}\varphi_{j}Z_{i-j}, \qquad i \in \mathbb{Z},$$
with coefficients satisfying
$$ \sum_{j=-\infty}^{\infty}|\varphi_{j}|^{\delta} < \infty \quad \textrm{for some} \ 0 < \delta < \alpha,\,\delta \leqslant 1.$$
 Astrauskas~\cite{As83} and Davis and Resnick~\cite{DaRe85} showed that the normalized sums of $X_i$'s under these conditions
converge in distribution to a stable random variable.
A natural generalization of this result would be a functional limit theorem for the partial sum process of $X_{i}$'s  with respect to some natural topology on $D[0,1]$. In other words, it is interesting to show
\be\label{eq:AvTaintr}
\frac{1}{a_{n}}\sum_{i=1}^{\floor{n\,\cdot}} (X_{i} - c_{n}) \dto \bigg( \sum_{j=-\infty}^{\infty}\varphi_{j} \bigg) V(\,\cdot\,),
\ee
in $D[0,1]$, where $V(\,\cdot\,)$ is an $\alpha$--stable L\'{e}vy process and $c_{n}$ are appropriate centering constants, $D[0,1]$
being the space of right continuous functions on $[0,1]$ with left limits.

 If $X_{i}$ is a finite order moving average with at least two nonzero coefficients, then the convergence in (\ref{eq:AvTaintr}) cannot hold in the $J_{1}$ sense, see Avram and Taqqu~\cite{AvTa92} for instance.
  However, if all coefficients $\varphi_{i}$ are nonnegative, then the convergence in (\ref{eq:AvTaintr}) holds in the $M_{1}$ topology
  according to
 Avram and Taqqu~\cite{AvTa92}, see their
 Theorem 2 (see also Corollary 1 in Tyran-Kami\'{n}ska~\cite{Ty10b}).

In the same article, Avram and Taqqu  made the following conjecture: if $\varphi_{j}=0$ for $j < 0$, $\varphi_{0}, \varphi_{1}, \ldots \in \mathbb{R}$ and if for every $K$,
$$ 0 \leqslant \sum_{j=0}^{K}\varphi_{j} \bigg/ \sum_{j=0}^{\infty}\varphi_{j} \leqslant 1,$$
then (\ref{eq:AvTaintr}) holds in  $M_{2}$ topology.
This topology, is again due to Skorohod (for an extensive discussion of
 topologies on $D[0,1]$ we refer to Whitt~\cite{Whitt02}).
As our main result  we
give a proof of  Avram and Taqqu's conjecture.
In order to do so, we first recall the precise definition of the $M_{2}$ topology.
We proceed by proving the conjecture for the finite order moving average processes in Section~\ref{S:FiniteMA}, and then finally in Section~\ref{S:InfiniteMA}, we
extend this to infinite order moving average processes.

The $M_{2}$ topology on $D[0, 1]$ is defined using completed graphs. For $x \in D[0,1]$ the \emph{completed graph} of $x$ is the set
\[
  \Gamma_{x}
  = \{ (t,z) \in [0,1] \times \mathbb{R} : z= \lambda x(t-) + (1-\lambda)x(t) \ \text{for some}\ \lambda \in [0,1] \},
\]
where $x(t-)$ is the left limit of $x$ at $t$. Besides the points of the graph $ \{ (t,x(t)) : t \in [0,1] \}$, the completed graph of $x$ also contains the vertical line segments joining $(t,x(t))$ and $(t,x(t-))$ for all discontinuity points $t$ of $x$.
An $M_{2}$ \emph{parametric representation} of the completed graph $\Gamma_{x}$ is a continuous function $(r,u)$ mapping $[0,1]$ onto $\Gamma_{x}$ such that $r$ is nondecreasing, with $r$ being the time component and $u$ being the spatial component.  Let $\Pi_{s,2}(x)$ denote the set of $M_{2}$ parametric representations of the graph $\Gamma_{x}$. For $x_{1},x_{2} \in D[0,1]$ define
\[
  d_{s,2}(x_{1},x_{2})
  = \inf \{ \|r_{1}-r_{2}\|_{[0,1]} \vee \|u_{1}-u_{2}\|_{[0,1]} : (r_{i},u_{i}) \in \Pi_{s,2}(x_{i}), i=1,2 \},
\]
where $\|x\|_{[0,1]} = \sup \{ |x(t)| : t \in [0,1] \}$ and $a \vee b = \max\{a,b\}$. Now we say that $x_{n} \to x$ in $D[0,1]$ for a sequence $(x_{n})$ in the Skorohod $M_{2}$ topology if $d_{s,2}(x_{n}, x) \to 0$ as $n \to \infty$. The $M_{2}$ topology is weaker than the more frequently used $M_{1}$ and $J_{1}$ topologies which are also due to Skorohod. The $M_{2}$ topology can be generated using the Hausdorff metric on the spaces of graphs. For $x_{1},x_{2} \in D[0,1]$ define
$$ d_{M_{2}}(x_{1}, x_{2}) = \bigg(\sup_{a \in \Gamma_{x_{1}}} \inf_{b \in \Gamma_{x_{2}}} d(a,b) \bigg) \vee \bigg(\sup_{a \in \Gamma_{x_{2}}} \inf_{b \in \Gamma_{x_{1}}} d(a,b) \bigg),$$
where $d$ is the metric on $\mathbb{R}^{2}$ defined by $d((x_{1},y_{1}),(x_{2},y_{2}))=|x_{1}-x_{2}| \vee |y_{1}-y_{2}|$ for $(x_{i},y_{i}) \in \mathbb{R}^{2},\,i=1,2$.
The metric $d_{M_{2}}$ induces the $M_{2}$ topology.

\section{Finite order MA processes}
\label{S:FiniteMA}

 Let $\varphi_{0}, \varphi_{1}, \ldots , \varphi_{q}$ (for some fixed $q \in \mathbb{N}$) be real numbers satisfying
\be\label{eq:FiniteMAcond}
0 \le \sum_{i=0}^{s}\varphi_{i} \Bigg/ \sum_{i=0}^{q}\varphi_{i} \le 1, \qquad \textrm{for every} \ s=0, 1, \ldots, q.
\ee
Put $ \Phi= \Phi(q) = \sum_{i=0}^{q}\varphi_{i}$. Without loss of generality assume $\Phi > 0$.
The case $\Phi<0$ is completely equivalent if we multiply the
noise sequence $(Z_i)$ by minus 1, and is therefore omitted.
Observe that condition (\ref{eq:FiniteMAcond}) implies
\be
 \nonumber \sum_{i=0}^{s}\varphi_{i} \ge 0 \quad \textrm{and} \quad \sum_{i=s}^{q}\varphi_{i} \ge 0, \qquad \textrm{for every} \ s=0,1,\ldots,q.
\ee

Let $(X_{t})$ be a moving average process defined by
$$ X_{t} = \sum_{i=0}^{q}\varphi_{i}Z_{t-i}, \qquad t \in \mathbb{Z}.$$
Define further the corresponding partial sum process
\be\label{eq:defVn}
V_{n}(t) = \frac{1}{a_{n}} \Bigg( \sum_{i=1}^{\floor {nt}}X_{i} - \floor {nt}b_{n}\Bigg), \qquad t \in [0,1],
\ee
where
$$ b_{n} = \left\{ \begin{array}{cc}
                                   0, & \quad \alpha \in (0,1] \\
                                   \Phi\,\mathrm{E}(Z_{1}), & \quad \alpha \in (1,2)
                                 \end{array}\right..$$

\begin{thm}\label{t:FinMA}
Let $(Z_{i})_{i \in \mathbb{Z}}$ be an i.i.d. sequence of regularly varying random variables with index $\alpha \in (0,2)$. When $\alpha=1$, suppose further that $Z_{1}$ is symmetric.
Assume real numbers $\varphi_{0}, \varphi_{1}, \ldots , \varphi_{q}$ satisfy (\ref{eq:FiniteMAcond}). Then
$$ V_{n}(\,\cdot\,) \dto \Phi V(\,\cdot\,), \qquad n \to \infty,$$
in $D[0,1]$ endowed with the $M_{2}$ topology, where $V$ is an $\alpha$--stable L\'{e}vy process.
\end{thm}

\begin{rem}\label{r:chartriple}
The characteristic L\'evy
 triple of the limiting process $V$ in the theorem is of the form $(0,\mu,b)$, with $\mu$ as in $(\ref{eq:mu})$ and
$$ b = \left\{ \begin{array}{cc}
                                   0, & \quad \alpha = 1\\[0.4em]
                                   (p-r)\frac{\alpha}{1-\alpha}, & \quad \alpha \in (0,1) \cup (1,2)
                                 \end{array}\right..$$
\end{rem}

In the proof of the theorem we are going to use the following
simple lemma.

\begin{lem}\label{l:first}
\begin{itemize}
  \item[(i)] For $k < q$ it holds
  \begin{eqnarray*}
  \sum_{i=1}^{k}\frac{\Phi\,Z_{i}}{a_{n}} - \sum_{i=1}^{k}\frac{X_{i}}{a_{n}} & = & \sum_{u=0}^{k-1}\frac{Z_{k-u}}{a_{n}} \sum_{s=u+1}^{q}\varphi_{s} - \sum_{u=k-q}^{q-1}\frac{Z_{-u}}{a_{n}} \sum_{s=u+1}^{q}\varphi_{s}\\[0.6em]
  & & - \sum_{u=0}^{q-k-1}\frac{Z_{-u}}{a_{n}} \sum_{s=u+1}^{u+k}\varphi_{s}.
  \end{eqnarray*}
  \item[(ii)] For $k \ge q$ it holds
  \begin{eqnarray*}
  \sum_{i=1}^{k}\frac{\Phi\,Z_{i}}{a_{n}} - \sum_{i=1}^{k}\frac{X_{i}}{a_{n}} & = & \sum_{u=0}^{q-1}\frac{Z_{k-u}}{a_{n}} \sum_{s=u+1}^{q}\varphi_{s} - \sum_{u=0}^{q-1}\frac{Z_{-u}}{a_{n}} \sum_{s=u+1}^{q}\varphi_{s}\\[0.7em]
  & =: & H_{n}(k) - G_{n}.
  \end{eqnarray*}
  \item[(iii)] For $q \le k \le n-q$ it holds
  \begin{eqnarray*}
  \sum_{i=1}^{k}\frac{\Phi\,Z_{i}}{a_{n}} - \sum_{i=1}^{k+q}\frac{X_{i}}{a_{n}} & = & - \sum_{u=0}^{q-1}\frac{Z_{-u}}{a_{n}} \sum_{s=u+1}^{q}\varphi_{s} - \sum_{u=1}^{q}\frac{Z_{k+u}}{a_{n}} \sum_{s=0}^{q-u}\varphi_{s}\\[0.7em]
  & =: & -G_{n} - T_{n}(k).
  \end{eqnarray*}
\end{itemize}
\end{lem}

\begin{pot1} We prove only (i), since the other two statements can be proven similarly. Note first that for every $k \in \mathbb{N}$ it holds that
\begin{eqnarray}\label{eq:lemmas1}
  \nonumber \sum_{i=1}^{k} \sum_{j=0}^{q} \varphi_{j}Z_{i-j} &=& \sum_{l=1-q}^{k} Z_{l} \sum_{i=1 \vee l}^{k \wedge (q+l)}\varphi_{i-l} \\
   &=& \sum_{l=1-q}^{0} Z_{l} \sum_{s=1-l}^{(k-l) \wedge q}\varphi_{s} + \sum_{l=1}^{k} Z_{l} \sum_{s=0}^{(k-l) \wedge q}\varphi_{s}
\end{eqnarray}
Since $k < q$, by (\ref{eq:lemmas1}) we have (recall $\Phi = \sum_{i=0}^{q}\varphi_{i}$)
\begin{eqnarray*}
  \sum_{i=1}^{k}\Phi\,Z_{i} - \sum_{i=1}^{k}X_{i} &  &  \\
   & \hspace*{-10em} = & \hspace*{-5em} \ \sum_{i=1}^{k}\Phi\,Z_{i} - \sum_{l=1-q}^{k-q}Z_{l} \sum_{s=1-l}^{q}\varphi_{s} - \sum_{l=k-q+1}^{0}Z_{l} \sum_{s=1-l}^{k-l}\varphi_{s} - \sum_{l=1}^{k} Z_{l} \sum_{s=0}^{k-l}\varphi_{s}\\
    & \hspace*{-10em} = & \hspace*{-5em}  \ \sum_{l=1}^{k}Z_{l} \sum_{s=k-l+1}^{q}\varphi_{s} - \sum_{l=1-q}^{k-q}Z_{l} \sum_{s=1-l}^{q}\varphi_{s} - \sum_{l=k-q+1}^{0}Z_{l} \sum_{s=1-l}^{k-l}\varphi_{s}.
\end{eqnarray*}
Now, we use the change of variables ($u=k-l$ for the first term on the right hand side in the last equation, and $u=-l$ for the second and third term) and rearrange some sums to arrive at
\begin{eqnarray*}
  \sum_{i=1}^{k}\Phi\,Z_{i} - \sum_{i=1}^{k}X_{i} & = & \sum_{u=0}^{k-1}Z_{k-u} \sum_{s=u+1}^{q}\varphi_{s} - \sum_{u=k-q}^{q-1}Z_{-u} \sum_{s=u+1}^{q}\varphi_{s}\\[0.6em]
  & & - \sum_{u=0}^{q-k-1}Z_{-u} \sum_{s=u+1}^{u+k}\varphi_{s}.
  \end{eqnarray*}\qed
\end{pot1}

\begin{rem}
Note that random variables $H_{n}(k)$ and $T_{n}(k)$ are independent.
\end{rem}

\begin{pot2} Case $\alpha \in (0,1]$.
Since the random variables $Z_{i}$ are i.i.d. and regularly varying,
Theorem 7.1 and Corollary 7.1 in Resnick~\cite{Resnick07}
 and Karamata's theorem immediately yield
$V_{n}^{Z}(\,\cdot\,) \dto \Phi\,V(\,\cdot\,)$, as $n \to \infty$, in $(D[0,1], d_{J_{1}})$, where
$$ V_{n}^{Z}(t) := \sum_{i=1}^{\floor {nt}}\frac{\Phi\,Z_{i}}{a_{n}}, \qquad t \in [0,1]$$
and
 $V$ is an $\alpha$--stable L\'{e}vy process
 with characteristic triple $(0,\mu,0)$ if $\alpha=1$ and $(0,\mu,(p-r)\alpha/(1-\alpha))$ if $\alpha \in (0,1)$
  with $p$ and $r$ as in (\ref{eq:pq}).

Using the fact that $J_{1}$ convergence implies $M_{2}$ convergence, we obtain
\be
V_{n}^{Z}(\,\cdot\,) \dto \Phi\,V(\,\cdot\,), \qquad n \to \infty,
\ee
in $(D[0,1], d_{M_{2}})$ as well. If one can show that for every $\epsilon >0$
$$ \lim_{n \to \infty}\Pr[d_{M_{2}}(V_{n}^{Z}, V_{n})> \epsilon]=0,$$
an application of  Slutsky's theorem (see for instance Theorem 3.4 in Resnick~\cite{Resnick07}),
will imply
$V_{n}(\,\cdot\,) \dto \Phi\,V(\,\cdot\,)$, as $n \to \infty$, in $(D[0,1], d_{M_{2}})$.

Fix $\epsilon >0$ and let $n \in \mathbb{N}$ be large enough, i.e.
 $n > \max\{2q, 4q/\epsilon\}$.
Then by the definition of the metric $d_{M_{2}}$, we have
\begin{eqnarray*}
  d_{M_{2}}(V_{n}^{Z},V_{n}) &=& \bigg(\sup_{a \in \Gamma_{V_{n}^{Z}}} \inf_{b \in \Gamma_{V_{n}}} d(a,b) \bigg) \vee \bigg(\sup_{a \in \Gamma_{V_{n}}} \inf_{b \in \Gamma_{V_{n}^{Z}}} d(a,b) \bigg) \\[0.4em]
   &= :& Y_{n} \vee T_{n}.
\end{eqnarray*}
Hence
\be\label{eq:AB}
\Pr [d_{M_{2}}(V_{n}^{Z}, V_{n})> \epsilon ] \leqslant \Pr(Y_{n}>\epsilon) + \Pr(T_{n}>\epsilon)
\ee
Now, we estimate the first term on the right hand side of (\ref{eq:AB}).
By the definition of $Y_{n}$, the Hausdorff metric and the choice of number $n$,
we see that
\begin{eqnarray}\label{eq:Yn}
  \nonumber\{Y_{n} > \epsilon\} & \subseteq & \{\exists\,a \in \Gamma_{V_{n}^{Z}} \ \textrm{such that} \ d(a,b) > \epsilon \ \textrm{for every} \ b \in \Gamma_{V_{n}} \} \\[0.6em]
  \nonumber & \subseteq & \{\exists\,k \in \{1,\ldots,q-1\} \ \textrm{such that} \ | V_{n}^{Z}(k/n) - V_{n}(k/n)| > \epsilon \}\\[0.6em]
  \nonumber & & \cup \ \{\exists\,k \in \{q,\ldots,n-q\} \ \textrm{such that} \ | V_{n}^{Z}(k/n) - V_{n}(k/n)| > \epsilon\\[0.6em]
  \nonumber & & \hspace*{1.5em} \textrm{and} \  | V_{n}^{Z}(k/n) - V_{n}((k+q)/n)| > \epsilon \}\\[0.6em]
  \nonumber & & \cup \ \{\exists\,k \in \{n-q+1,\ldots,n\} \ \textrm{such that} \ | V_{n}^{Z}(k/n) - V_{n}(k/n)| > \epsilon \}\\[0.6em]
  & =: & A^{Y}_{n} \cup B^{Y}_{n} \cup C^{Y}_{n}.
\end{eqnarray}
By Lemma~\ref{l:first} (i) and stationarity we obtain
\begin{eqnarray}\label{eq:Bnlemmafirst}
  \nonumber \Pr (A^{Y}_{n}) & \leqslant &  \sum_{k=1}^{q-1} \Pr \Big( \Big| \sum_{i=1}^{k}\frac{\Phi\,Z_{i}}{a_{n}} - \sum_{i=1}^{k}\frac{X_{i}}{a_{n}} \Big| >\epsilon \Big) \\[0.6em]
    \nonumber  & \leqslant & \sum_{k=1}^{q-1} \bigg[ \Pr \Big( \sum_{u=0}^{k-1}\frac{|Z_{k-u}|}{a_{n}}
         \sum_{s=u+1}^{q}|\varphi_{s}| > \frac{\epsilon}{3} \Big) + \Pr \Big( \sum_{u=k-q}^{q-1}\frac{|Z_{-u}|}{a_{n}} \sum_{s=u+1}^{q}|\varphi_{s}| > \frac{\epsilon}{3} \Big)\\[0.6em]
  \nonumber  & & \hspace*{2em} + \Pr \Big( \sum_{u=0}^{q-k-1}\frac{|Z_{-u}|}{a_{n}} \sum_{s=u+1}^{u+k}|\varphi_{s}| > \frac{\epsilon}{3} \Big) \bigg]\\[0.6em]
      & \leqslant &  3(q-1)(2q-1) \Pr \Big( \frac{|Z_{0}|}{a_{n}} > \frac{\epsilon}{3(2q-1)\theta} \Big),
  \end{eqnarray}
   where $\theta = \sum_{s=0}^{q}|\varphi_{s}| >0$. Hence, by
 regular variation property
we observe
 \be\label{eq:setBn1}
 \lim_{n \to \infty} \Pr (A^{Y}_{n}) = 0.
 \ee
Next, using Lemma~\ref{l:first} (ii), (iii) together with
stationarity and the fact that $H_{n}(k)$ and $T_{n}(k)$ are independent,
we obtain
\begin{eqnarray*}
  \Pr (B^{Y}_{n}) & = &  \Pr \Big( \exists\,k \in \{q,\ldots,n-q\} \ \textrm{such that} \ |H_{n}(k)-G_{n}| > \epsilon \\[0.6em]
  & & \hspace*{2em} \textrm{and} \ |-G_{n}-T_{n}(k)| > \epsilon \Big)\\[0.6em]
   & \leqslant & \Pr \Big( |G_{n}|> \frac{\epsilon}{2} \Big) + \sum_{k=q}^{n-q} \Pr \Big( |H_{n}(k)| >
       \frac{\epsilon}{2} \ \textrm{and} \ |T_{n}(k)| > \frac{\epsilon}{2} \Big)\\[0.6em]
  &= & \Pr \Big( |G_{n}|> \frac{\epsilon}{2} \Big) + \sum_{k=q}^{n-q} \Pr \Big( |H_{n}(k)| > \frac{\epsilon}{2} \Big)
           \Pr \Big(|T_{n}(k)| > \frac{\epsilon}{2} \Big)\\[0.6em]
  & \leqslant & \Pr \Big( |G_{n}|> \frac{\epsilon}{2} \Big) + n \Pr \Big( |H_{n}(0)| > \frac{\epsilon}{2} \Big) \Pr
            \Big(|T_{n}(0)| > \frac{\epsilon}{2} \Big)\\[0.6em]
   &\leqslant &  q \Pr \Big( \frac{|Z_{0}|}{a_{n}} > \frac{\epsilon}{2q\theta} \Big) + \frac{q^{2}}{n} \Big[ n \Pr
            \Big( \frac{|Z_{0}|}{a_{n}} > \frac{\epsilon}{2q\theta} \Big) \Big]^{2},
  \end{eqnarray*}
whence we conclude
 \be\label{eq:setBn3}
 \lim_{n \to \infty} \Pr (B^{Y}_{n}) = 0.
 \ee
In a similar manner as in (\ref{eq:Bnlemmafirst}), but using (ii) from Lemma~\ref{l:first} instead of (i) we get
  \be\label{eq:setBn4}
   \lim_{n \to \infty} \Pr (C^{Y}_{n})=0.
  \ee
  From relations (\ref{eq:Yn}), (\ref{eq:setBn1}), (\ref{eq:setBn3}) and (\ref{eq:setBn4}) we obtain
  \be\label{eq:Ynend}
  \lim_{n \to \infty} \Pr(Y_{n} > \epsilon ) =0.
  \ee

 It remains to estimate the second term on the right hand side of (\ref{eq:AB}). From the definition of $T_{n}$, the Hausdorff metric and the number $n$ it follows
\begin{eqnarray}\label{eq:Zn}
  \nonumber\{T_{n} > \epsilon\} & \subseteq & \{\exists\,a \in \Gamma_{V_{n}} \ \textrm{such that} \ d(a,b) > \epsilon \ \textrm{for every} \ b \in \Gamma_{V_{n}^{Z}} \} \\[0.6em]
  \nonumber & \subseteq & \{\exists\,k \in \{1,\ldots,2q-1\} \
\ \textrm{such that} \ | V_{n}(k/n) - V_{n}^{Z}(k/n)| > \epsilon \}\\[0.6em]
  \nonumber & & \cup \ \Big\{\exists\,k \in \{2q,\ldots,n\} \
   \textrm{such that} \
  d((k/n, V_{n}(k/n)), \Gamma_{V_{n}^{Z}})> \frac{\epsilon}{2}
    \Big\}\\[0.6em]
  & =: & A^{T}_{n} \cup B^{T}_{n}.
\end{eqnarray}
Using Lemma~\ref{l:first} (i) and (ii), one could similarly as before for set $A^{Y}_{n}$ obtain
\be\label{eq:Cnfirst}
\lim_{n \to \infty} \Pr( A^{T}_{n})=0.
\ee
To bound $\Pr(B^{T}_{n})$ we need a new argument.
 For each $k \geqslant 2q$, set $V^{Z,\min}_k = \min\{V^Z_n((k-q)/n), V^Z_n(k/n) \}$
and $V^{Z,\max}_k = \max\{V^Z_n((k-q)/n), V^Z_n(k/n) \}$.
Since the completed graph $\Gamma_{V_{n}^{Z}}$ is connected,
if
$$
 V_n(k/n) \in (V^{Z,\min}_k - \frac{\epsilon}{4},V^{Z,\max}_k + \frac{\epsilon}{4} )
$$
then
$d((k/n, V_{n}(k/n)), \Gamma_{V_{n}^{Z}})<\epsilon/2$ for all $n$ large enough so that $q/n < \epsilon/4$.
Therefore, $P(B_n^T)$ is bounded by
\begin{eqnarray*}
\lefteqn{P \left(
\exists\,k \in \{2q,\ldots,n\} \ \textrm{such that} \
\sum_{i=1}^k \frac{X_i}{a_n} > V^{Z,\max}_k + \frac{\epsilon}{4} \right)}\\
&+&
P \left(
\exists\,k \in \{2q,\ldots,n\} \ \textrm{such that} \
\sum_{i=1}^k \frac{X_i}{a_n} < V^{Z,\min}_k - \frac{\epsilon}{4} \right)
\end{eqnarray*}
In the sequel we consider only the first of these two probabilities,
since the other one can be handled in a similar manner.
Note, that the first probability
using Lemma~\ref{l:first},  can be bounded  by
\begin{eqnarray*}
\lefteqn{P \left(
\exists\,k \in \{2q,\ldots,n\} \ \textrm{such that} \
G_n - H_n(k) > \frac{\epsilon}{4} \ \mbox{ and }\
G_n + T_n(k-q) > \frac{\epsilon}{4}   \right)} \\
& \leq &
P\left(G_n  > \frac{\epsilon}{8}\right)\\
&&  +
P \left(
\exists\,k \in \{2q,\ldots,n\} \ \textrm{such that} \
H_n(k) < - \frac{\epsilon}{8} \ \mbox{ and }\
T_n(k-q) > \frac{\epsilon}{8}   \right)
\end{eqnarray*}
As before $P(G_n  > {\epsilon}/{8})\to 0$
 as $n \to \infty$. For the second term,
note,
$$
 H_n(k)= \sum_{u=0}^{q-1}\frac{Z_{k-u}}{a_{n}} \sum_{s=u+1}^{q}\varphi_{s}
\ \mbox{ and } T_n(k-q) = \sum_{u=0}^{q-1}\frac{Z_{k-u}}{a_{n}}
 \sum_{s=0}^{u}\varphi_{s}\,.
$$
Hence, that term is bounded by
\[
n P\left(
\sum_{u=0}^{q-1}\frac{Z_{-u}}{a_{n}} \sum_{s=u+1}^{q}\varphi_{s}
< - \frac{\epsilon}{8}
\ \mbox{ and }
\sum_{u=0}^{q-1}\frac{Z_{-u}}{a_{n}}
 \sum_{s=0}^{u}\varphi_{s}
> \frac{\epsilon}{8}
 \right)
\]
where we used the stationarity of the sequence $(Z_i)$ .
Observe now that the sums
$\sum_{s=0}^{u}\varphi_{s}$
and $ \sum_{s=u+1}^{q}\varphi_{s}$ are both nonnegative and
bounded by $\Phi =  \sum_{s=0}^{q}\varphi_{s}$, see \eqref{eq:FiniteMAcond}.
Therefore, the last expression above is bounded by
\begin{eqnarray*}
& &{n P\left(
\exists\, i,j \in \{0,\ldots,q-1\},\,i \neq j \ \textrm{such that} \
\Phi \frac{Z_{-i }}{a_{n}}< - \frac{\epsilon}{8q}
 \mbox{ and } \Phi \frac{Z_{-j }}{a_{n}} >\frac{\epsilon}{8q}  \right)} \\
& \leq &  n {q \choose 2}
 P\left(
 \frac{|Z_{0}|}{a_{n}} >  \frac{\epsilon}{8 q \Phi}
 \right)^2,
\end{eqnarray*}
which clearly tends to 0 as $n\toi$, by the regular variation property of
the random variables $Z_i$.
 Note that the case $i=j$ above is not possible since then we would have $Z_{-i}<0$ and $Z_{-i}>0$.

Together with relations (\ref{eq:Zn}) and (\ref{eq:Cnfirst}) this implies
\be\label{eq:Tnend}
\lim_{n \to \infty} \Pr(T_{n}>\epsilon)=0.
\ee
Now from (\ref{eq:AB}), (\ref{eq:Ynend}) and (\ref{eq:Tnend}) we obtain
\be
\lim_{n \to \infty} \Pr [d_{M_{2}}(V_{n}^{Z}, V_{n})> \epsilon ]=0,
\ee
and finally we conclude that $V_{n}(\,\cdot\,) \dto \Phi\,V(\,\cdot\,)$, as $n \to \infty$, in $(D[0,1], d_{M_{2}})$.\\

Case $\alpha \in (1,2)$. In this case $\mathrm{E}(Z_{1}) < \infty$. Define
$$ Z_{i}' = Z_{i} - \mathrm{E}(Z_{1}), \qquad i \in \mathbb{Z}.$$
Then $\mathrm{E}(Z_{i}')=0$ and $(Z_{i}')_{i}$ is an i.i.d. sequence of regularly varying random variables with index $\alpha$. Then it is known that, as $n \to \infty$, the stochastic process
\begin{equation*}\label{eq:AOPthm}
 W_{n}(t) := \sum_{i=1}^{\floor {nt} }\frac{Z_{i}}{a_{n}} - \floor{nt} \mathrm{E} \Big(\frac{Z_{1}}{a_{n}} 1_{\{|Z_{1}| \leqslant a_{n}\}} \Big), \qquad t \in [0,1],
\end{equation*}
converges in distribution in $(D[0,1], d_{M_{1}})$ to an $\alpha$--stable L\'{e}vy process with characteristic triple $(0,\mu,0)$ (cf. Theorem 3.4 in Basrak et al.~\cite{BKS}). By Karamata's theorem, as $n \to \infty$,
$$ n\,\mathrm{E} \Big( \frac{Z_{1}}{a_{n}} 1_{\{ |Z_{1}| > a_{n} \}} \Big) \to (p-r)\frac{\alpha}{\alpha-1},$$
with $p$ and $r$ as in (\ref{eq:pq}). Thus, as $n \to \infty$,
$$ \floor{n\,\cdot} \mathrm{E} \Big( \frac{Z_{1}}{a_{n}} 1_{\{ |Z_{1}| > a_{n} \}} \Big) \to (\,\cdot\,) (p-r)\frac{\alpha}{\alpha-1}$$
in $(D[0,1], d_{M_{1}})$.
Since the latter function is continuous, an application of Corollary 12.7.1 in Whitt~\cite{Whitt02} (which gives a sufficient condition for addition to be continuous) and the continuous mapping theorem give that the following stochastic process
$$ \sum_{i=1}^{\floor {nt} }\frac{Z_{i}'}{a_{n}} = \sum_{i=1}^{\floor {nt} }\frac{Z_{i}}{a_{n}} - \floor{nt} \mathrm{E} \Big(\frac{Z_{1}}{a_{n}} \Big)
= W_{n}(t) - \floor{nt} \mathrm{E} \Big( \frac{Z_{1}}{a_{n}} 1_{\{ |Z_{1}| > a_{n} \}} \Big), \quad t \in [0,1],$$
converges in distribution in $(D[0,1], d_{M_{1}})$ to an $\alpha$--stable L\'{e}vy process $V$ with characteristic triple $(0,\mu,(p-r)\alpha/(1-\alpha))$. Define now
$$ X_{i}' = \sum_{j=0}^{q}\varphi_{j}Z_{i-j}', \qquad i \in \mathbb{Z},$$
and
$$V_{n}'(t) = \sum_{i=1}^{\floor{nt}}\frac{X_{i}'}{a_{n}} \quad \textrm{and} \quad V_{n}'^{Z}(t) = \sum_{i=1}^{\floor{nt}}\frac{\Phi\,Z_{i}'}{a_{n}}, \quad t \in [0,1].$$
Now we can repeat all arguments used in the case $\alpha \in (0,1]$ to obtain
\be
V_{n}'^{Z}(\,\cdot\,) \dto \Phi\,V(\,\cdot\,), \qquad n \to \infty,
\ee
in $(D[0,1], d_{M_{2}})$, and
$$ \lim_{n \to \infty}\Pr[d_{M_{2}}(V_{n}'^{Z}, V_{n}')> \epsilon]=0, \qquad \epsilon >0.$$
Note that $V_{n} = V_{n}'$ and therefore $V_{n}(\,\cdot\,) \dto \Phi\,V(\,\cdot\,)$, as $n \to \infty$, in $(D[0,1], d_{M_{2}})$. This concludes the proof.\qed
\end{pot2}

\section{Infinite order MA processes}
\label{S:InfiniteMA}

Let $\{\varphi_{i}, i=0,1,2,\ldots\}$ be a sequence of real numbers satisfying
\be\label{eq:convcond}
 \sum_{i=0}^{\infty}|\varphi_{i}|^{\delta} < \infty
\ee
for some $0 < \delta < \min\{1,\alpha\}$, and
\be\label{eq:InfiniteMAcond}
0 \le \sum_{i=0}^{s}\varphi_{i} \Bigg/ \sum_{i=0}^{\infty}\varphi_{i} \le 1, \qquad \textrm{for every} \ s=0, 1, 2 \ldots.
\ee
Let $ \Phi= \Phi(\infty) = \sum_{i=0}^{\infty}\varphi_{i}$. Condition (\ref{eq:convcond}) implies $\Phi$ is finite. Without loss of generality assume $\Phi > 0$
 (as before, the case $\Phi <0$ can be handled similarly).

Let $(X_{t})$ be a moving average process defined by
$$ X_{t} = \sum_{i=0}^{\infty}\varphi_{i}Z_{t-i}, \qquad t \in \mathbb{Z}.$$
Condition (\ref{eq:convcond}) ensures that $X_{t}$ converges in $L^{\delta}$ and a.s.
Define further the corresponding partial sum stochastic process $V_{n}$ as in (\ref{eq:defVn}).

\begin{thm}\label{t:InfMA}
Let $(Z_{i})_{i \in \mathbb{Z}}$ be an i.i.d. sequence of regularly varying random variables with index $\alpha \in (0,2)$. When $\alpha=1$, suppose further that $Z_{1}$ is symmetric. Let $\{\varphi_{i}, i=0,1,2,\ldots\}$ be a sequence of real numbers satisfying (\ref{eq:convcond}) and (\ref{eq:InfiniteMAcond}). Then
$$ V_{n}(\,\cdot\,) \dto \Phi V(\,\cdot\,), \qquad n \to \infty,$$
in $D[0,1]$ endowed with the $M_{2}$ topology, where $V$ is an $\alpha$--stable L\'{e}vy process.
\end{thm}

\begin{rem}
The characteristic triple of the limiting process $V$ in Theorem~\ref{t:InfMA} is of the same form as in Remark~\ref{r:chartriple}.
\end{rem}

\begin{pot3} Case $\alpha \in (0,1]$. Fix $q \in \mathbb{N}$ and define
$$ X_{i}^{q} = \sum_{j=0}^{q-1}\varphi_{j}Z_{i-j} + \varphi'_{q} Z_{i-q}, \qquad i \in \mathbb{Z},$$
where $\varphi'_{q}= \sum_{i=q}^{\infty}\varphi_{i}$,
and
$$ V_{n, q}(t) = \sum_{i=1}^{\floor{nt}} \frac{X_{i}^{q}}{a_{n}}, \qquad t \in [0,1].$$
 Since the coefficients $\varphi_{0}, \ldots, \varphi_{q-1}, \varphi'_{q}$ satisfy condition (\ref{eq:FiniteMAcond}), an application of Theorem~\ref{t:FinMA} to a finite order moving average process $(X_{i}^{q})_{i}$ yields that
\be
V_{n, q}(\,\cdot\,) \dto \Phi\,V(\,\cdot\,), \qquad n \to \infty,
\ee
in $(D[0,1], d_{M_{2}})$. If we show that for every $\epsilon >0$
$$ \lim_{q \to \infty} \limsup_{n \to \infty}\Pr[d_{M_{2}}(V_{n, q}, V_{n})> \epsilon]=0,$$
then by a generalization of Slutsky's theorem (see for instance Theorem 3.5 in Resnick~\cite{Resnick07}) it will follow $V_{n}(\,\cdot\,) \dto \Phi\,V(\,\cdot\,)$, as $n \to \infty$, in $(D[0,1], d_{M_{2}})$. Since the Skorohod $M_{2}$ metric on $D[0,1]$ is bounded above by the uniform metric on $D[0,1]$, it suffices to show that
$$ \lim_{q \to \infty} \limsup_{n \to \infty}\Pr \bigg( \sup_{0 \leqslant t \leqslant 1}|V_{n, q}(t) - V_{n}(t)|> \epsilon \bigg)=0.$$
Recalling the definitions, we have
\begin{eqnarray*}
  \lim_{q \to \infty} \limsup_{n \to \infty}\Pr \bigg( \sup_{0 \leqslant t \leqslant 1}|V_{n, q}(t) - V_{n}(t)|> \epsilon \bigg) & &  \\[0.6em]
   & \hspace*{-26em} \leqslant & \hspace*{-13em} \lim_{q \to \infty} \limsup_{n \to \infty}\Pr \bigg( \sum_{i=1}^{n}\frac{|X_{i}^{q}-X_{i}|}{a_{n}} > \epsilon \bigg).
\end{eqnarray*}
Put $\varphi''_{q} = \varphi'_{q} - \varphi_{q} = \sum_{j=q+1}^{\infty}\varphi_{j}$ and observe
\begin{eqnarray*}
  \sum_{i=1}^{n}|X_{i}^{q}-X_{i}| & = & \sum_{i=1}^{n} \bigg| \sum_{j=0}^{q-1}\varphi_{j}Z_{i-j} + \varphi'_{q}Z_{i-q} - \sum_{j=0}^{\infty}\varphi_{j}Z_{i-j}\bigg| \\[0.6em]
  & = & \sum_{i=1}^{n} \bigg| \varphi''_{q} Z_{i-q} - \sum_{j=q+1}^{\infty}\varphi_{j}Z_{i-j}\bigg|\\[0.6em]
  & \leqslant & \sum_{i=1}^{n} \bigg[ |\varphi''_{q}|\,|Z_{i-q}| + \sum_{j=q+1}^{\infty}|\varphi_{j}|\,|Z_{i-j}| \bigg]\\[0.6em]
   & \leqslant & \bigg( 2 \sum_{j=q+1}^{\infty}|\varphi_{j}| \bigg) \sum_{i=1}^{n}|Z_{i-q}| + \sum_{i=-\infty}^{0}|Z_{i-q}| \sum_{j=1}^{n}|\varphi_{q-i+j}|\\[0.6em]
   & =: & O_{n}(q).
\end{eqnarray*}
Therefore we have to show
\be\label{eq:Oiqn1}
\lim_{q \to \infty} \limsup_{n \to \infty} \Pr \bigg(\frac{O_{n}(q)}{a_{n}} > \epsilon \bigg)=0.
\ee
By Lemma 2 in in Avram and Taqqu~\cite{AvTa92}, we obtain, for large $q$,
\begin{eqnarray*}
  \Pr \bigg(\frac{O_{n}(q)}{a_{n}} > \epsilon \bigg) & \leqslant & M \frac{\epsilon^{-(\alpha + \eta)}}{n} \bigg( n \Big( 2 \sum_{j=q+1}^{\infty}|\varphi_{j}| \Big)^{\alpha - \eta} + \sum_{i=-\infty}^{0} \Big( \sum_{j=1}^{n}|\varphi_{q-i+j}| \Big)^{\alpha - \eta} \bigg),
\end{eqnarray*}
where $\eta$ is some positive real number satisfying $\alpha - \eta > \delta$
 and $M$ is a positive constant.
Since $\alpha - \eta < \alpha \leqslant 1$, an application of the inequality $|\sum_{i=1}^{n}a_{i}|^{\gamma} \leqslant \sum_{i=1}^{n}|a_{i}|^{\gamma}$ with $a_{i}, \ldots, a_{n}$ real numbers and $\gamma \in (0,1]$, yields
$$\bigg( \sum_{j=1}^{n}|\varphi_{q-i+j}| \bigg)^{\alpha - \eta}  \leqslant \sum_{j=1}^{n}|\varphi_{q-i+j}|^{\alpha - \eta}.$$
Using this and the fact that every $|\varphi_{i}|^{\alpha - \eta}$, for $i=q+1, q+2, \ldots$, appears in the sum $\sum_{i=-\infty}^{0} \sum_{j=1}^{n}|\varphi_{q-i+j}|^{\alpha - \eta}$ at most $n$ times, we obtain
\begin{eqnarray}\label{eq:Oiqn2}
  \nonumber \Pr \bigg(\frac{O_{n}(q)}{a_{n}} > \epsilon \bigg) & \leqslant & M \frac{\epsilon^{-(\alpha + \eta)}}{n} \bigg( n \Big( 2 \sum_{j=q+1}^{\infty}|\varphi_{j}| \Big)^{\alpha - \eta} + n \sum_{j=q+1}^{\infty}|\varphi_{j}|^{\alpha - \eta} \bigg)\\[0.6em]
  & = & M \epsilon^{-(\alpha + \eta)} \bigg( \Big( 2 \sum_{j=q+1}^{\infty}|\varphi_{j}| \Big)^{\alpha - \eta} + \sum_{j=q+1}^{\infty}|\varphi_{j}|^{\alpha - \eta} \bigg).
\end{eqnarray}
 Since for large $q$ it holds that $|\varphi_{j}| \leqslant |\varphi_{j}|^{\delta}$ and $|\varphi_{j}|^{\alpha - \eta} \leqslant |\varphi_{j}|^{\delta}$ for all $j \geqslant q+1$, from condition (\ref{eq:convcond}) we immediately obtain, as $q \to \infty$,
$$ \sum_{j=q+1}^{\infty}|\varphi_{j}| \to 0 \quad \textrm{and} \quad \sum_{j=q+1}^{\infty}|\varphi_{j}|^{\alpha - \eta} \to 0.$$
Therefore from (\ref{eq:Oiqn2}) letting $q \to \infty$, follows (\ref{eq:Oiqn1}), which means that $V_{n}(\,\cdot\,) \dto \Phi\,V(\,\cdot\,)$, as $n \to \infty$, in $(D[0,1], d_{M_{2}})$.\\

Case $\alpha \in (1,2)$. Define
$ Z_{i}' = Z_{i} - \mathrm{E}(Z_{1}), \ i \in \mathbb{Z}$.
Fix $q \in \mathbb{N}$ and define
$$ X_{i}'^{q} = \sum_{j=0}^{q-1}\varphi_{j}Z'_{i-j} + \varphi'_{q} Z'_{i-q}, \qquad i \in \mathbb{Z},$$
and
$$ V_{n, q}'(t) = \sum_{i=1}^{\floor{nt}} \frac{X_{i}'^{q}}{a_{n}}, \qquad t \in [0,1].$$
Then
$$  V_{n, q}'(t) = \frac{1}{a_{n}} \bigg( \sum_{i=1}^{\floor{nt}}X_{i}^{q} - \floor{nt} \Phi \mathrm{E}(Z_{1}) \bigg).$$
 Since the coefficients $\varphi_{0}, \ldots, \varphi_{q-1}, \varphi'_{q}$ satisfy condition (\ref{eq:FiniteMAcond}), Theorem~\ref{t:FinMA}, applied to a finite order moving average process $(X_{i}'^{q})_{i}$ yields that
\be
V_{n, q}'(\,\cdot\,) \dto \Phi\,V(\,\cdot\,), \qquad n \to \infty,
\ee
in $(D[0,1], d_{M_{2}})$. In order to obtain $V_{n}(\,\cdot\,) \dto \Phi\,V(\,\cdot\,)$, as in the previous case, it remains to show that for every $\epsilon >0$
$$ \lim_{q \to \infty} \limsup_{n \to \infty}\Pr[d_{M_{2}}(V_{n, q}', V_{n})> \epsilon]=0,$$
i.e.
\begin{equation*}
\lim_{q \to \infty} \limsup_{n \to \infty} \Pr \bigg(\frac{O_{n}(q)}{a_{n}} > \epsilon \bigg)=0.
\end{equation*}
As before, by Lemma 2 in Avram and Taqqu~\cite{AvTa92}, for large $q$,
\begin{eqnarray*}
  \Pr \bigg(\frac{O_{n}(q)}{a_{n}} > \epsilon \bigg) & \leqslant & M \frac{\epsilon^{-(\alpha + \eta)}}{n} \bigg( n \Big( 2 \sum_{j=q+1}^{\infty}|\varphi_{j}| \Big)^{\alpha - \eta} + \sum_{i=-\infty}^{0} \Big( \sum_{j=1}^{n}|\varphi_{q-i+j}| \Big)^{\alpha - \eta} \bigg),
\end{eqnarray*}
where $\eta$ is some positive real number satisfying $\alpha - \eta >1$.
 Now using the inequality $|\sum_{i=1}^{n}a_{i}|^{\gamma} \leqslant \sum_{i=1}^{n}|a_{i}|$ with $a_{i}, \ldots, a_{n}$ real numbers such that $|a_{1} + \ldots + a_{n}|<1$  and $\gamma \in (1,2)$,
 similarly as before we obtain
\begin{eqnarray*}
  \nonumber \Pr \bigg(\frac{O_{n}(q)}{a_{n}} > \epsilon \bigg) & \leqslant & M \frac{\epsilon^{-(\alpha + \eta)}}{n} \bigg( n \Big( 2 \sum_{j=q+1}^{\infty}|\varphi_{j}| \Big)^{\alpha - \eta} + n \sum_{j=q+1}^{\infty}|\varphi_{j}| \bigg)\\[0.6em]
  & = & M \epsilon^{-(\alpha + \eta)} \bigg( \Big( 2 \sum_{j=q+1}^{\infty}|\varphi_{j}| \Big)^{\alpha - \eta} +  \sum_{j=q+1}^{\infty}|\varphi_{j}| \bigg).
\end{eqnarray*}
and again letting $\lim_{q \to \infty} \limsup_{n \to \infty}$, the desired result follows. This completes the proof. \qed
\end{pot3}

\section*{Acknowledgements}
Bojan Basrak's research was partially supported by the research grant MZOS nr.
037-0372790-2800 of the Croatian government.

\end{document}